\documentclass[a4paper,11pt,reqno,oneside]{article}

\usepackage{amsmath}
\usepackage{amsthm}
\usepackage{amssymb}

\usepackage{enumerate}
\usepackage{epic,eepic}

\usepackage[dvips]{graphicx}
\usepackage{calc}

\newcommand{\Z}{{\mathbb Z}}
\newcommand{\Q}{{\mathbb Q}}
\newcommand{\N}{{\mathbb N}}

\newcommand{\Acal}{\mathcal{A}}
\newcommand{\Mcal}{\mathcal{M}}
\newcommand{\Rcal}{\mathcal{R}}

\newcommand{\mat}[1]{\boldsymbol{#1}}
\newcommand{\coloneq}{\mathrel{\mathop:}=}
\newcommand{\eqcolon}{\mathrel{=\!\!\mathop:}}

\newcommand{\EtN}{\mathrm{E}(3,\N)}
\newcommand{\PhiT}{\Phi_{3\mathrm{iet}}}

\newtheorem{lem}{Lemma}
\newtheorem{thm}[lem]{Theorem}
\newtheorem{prop}[lem]{Proposition}
\newtheorem{coro}[lem]{Corollary}
\theoremstyle{definition}
\newtheorem{ex}[lem]{\emph{Example}}
\newtheorem*{defi}{\emph{Definition}}

\begin{document}

\title{\textbf{A note on 3iet preserving morphisms}}
\author{P.~Ambro\v{z},\quad E.~Pelantov\'a \\[5mm]
{\small Doppler Institute \& Department of Mathematics}\\
{\small FNSPE, Czech Technical University, Trojanova 13, 120 00 Praha 2, Czech Republic}\\
{\small E-mail: \texttt{petr.ambroz@fjfi.cvut.cz}, {\small\texttt{pelantova@km1.fjfi.cvut.cz}}}}
\date{}
\maketitle

\begin{abstract}
An infinite word, which is aperiodic and codes the orbit of a transformation of the exchange
of three intervals is called 3iet word. Such a word is thus a natural generalization of a
sturmian word to a word over 3-letter alphabet.
A morphism is said to be 3iet preserving if it maps any 3iet word to another 3iet word.
It is known that the monoid of morphisms preserving sturmian words is finitely generated.
On the contrary, in this note we prove that the monoid
of 3iet preserving morphisms is not finitely generated, that is,
there are infinitely many 3iet preserving morphisms, which cannot be written as a
non-trivial decomposition of other 3iet preserving morphisms.
\end{abstract}

\section{Introduction}

One of possible generalizations of sturmian words are the words coding orbits of $r$-interval
exchange transformations, called $r$iet words. These transformations have been introduced
by Katok and Stepin~\cite{katok-stepin-umn-22} and then studied by other
authors~\cite{boshernitzan-carroll-jam-72,keane-mz-141,rauzy-aa-34}.

A word coding an exchange of two intervals is either periodic or
sturmian. Therefore the combinatorial properties of infinite words
coding an exchange of two intervals are well described, in contrast
to the $r$iet words with $r\geq 3$, which were studied only a little. The only
result known for a general $r$ is the estimate on values of the factor
complexity, namely $\mathcal{C}(n)\leq (r-1)n+1$; other properties remain hidden
for general $r$.

An important step towards the characterization of words coding 3-interval
exchange transformations was accomplished by Ferenczi, Holton and Zamboni in the series of
papers~\cite{ferenczi-holton-zamboni-aif-51,ferenczi-holton-zamboni-jam-89,ferenczi-holton-zamboni-jam-93}.
Our paper is devoted to the study of morphisms, which preserve the set of 3iet words
(the reader is referred to Section~\ref{sec:prelim} for precise definitions).
Recall that morphisms preserving sturmian words were completely described
in~\cite{berstel-seebold-bbms-1,mignosi-seebold-jtnb-5,seebold-tcs-88}. Berstel, Mignosi and
S\'e\'ebold studied morphisms $\varphi$ such that $\varphi(u)$ is sturmian for each
sturmian word $u$. They showed that the monoid of such morphisms is finitely generated
by three morphisms.

The aim of this note is to show that the monoid of morphisms, which preserve 3iet words,
is not finitely generated; strictly speaking we show that there are infinitely many morphisms
$\psi$ such that
\begin{enumerate}
\item
  $\psi(u)$ is a 3iet word for all 3iet words $u$, i.e., $\psi$ is 3iet preserving
\item
  $\psi$ cannot be written as a non-trivial composition of other 3iet preserving morphisms.
\end{enumerate}

In the paper~\cite{wen-wen-crasp-318} it is showed that the monoid of morphisms preserving
sturmian words coincides with the set of positive invertible substitutions over a 2-letter
alphabet. Invertible substitutions over a 3-letter alphabet are characterized
in~\cite{tan-wen-zhang-aam-32}. Even though the set of invertible substitutions over
a 3-letter alphabet is not finitely generated, the monoid of their matrices is. It means that
the monoid of positive invertible substitutions over a 3-letter alphabet and
the monoid of 3iet preserving morphisms are different.

\section{Preliminaries}\label{sec:prelim}
\subsection{3iet words}

In this paper we deal with finite and infinite words over a finite alphabet
$\Acal = \{a_1,\ldots,a_k\}$. The set of all finite words over $\Acal$ is denoted by $\Acal^*$.
This set, equipped with the concatenation as a binary operation, is a free monoid having the
empty word as its identity.
The set of two-sided infinite words over an alphabet $\Acal$,
i.e., of two-sided infinite sequences of letters of $\Acal$, is
denoted by $\Acal^\Z$, its elements are words $u=(u_n)_{n\in\Z}$.

We concentrate on words coding an exchange of three intervals.

\begin{defi}
  Let $\alpha,\beta,\gamma$ be three positive real numbers. Denote
  $$
  \begin{array}{rcl}
  I_A&\coloneq&[0,\alpha)\\
  I_B&\coloneq& [\alpha,\alpha+\beta)\\
  I_C&\coloneq& [\alpha+\beta,\alpha+\beta+\gamma)
  \end{array}
  \quad\hbox{or }\quad
  \begin{array}{rcl}
  I_A&\coloneq& (0,\alpha]\\
  I_B&\coloneq& (\alpha,\alpha+\beta]\\
  I_C&\coloneq& (\alpha+\beta,\alpha+\beta+\gamma]
  \end{array}
  $$
  respectively, and $I:=I_A\cup I_B\cup I_C$.
  A mapping
  $T:I\rightarrow I$, given by
  \begin{equation}\label{eq:3iet}
  T(x) = \begin{cases}
  x + \beta+\gamma & \text{if $x\in I_A$,} \\
  x -\alpha+\gamma & \text{if $x\in I_B$,} \\
  x -\alpha-\beta & \text{if $x\in I_C$,}
  \end{cases}
  \end{equation}
  is called a \emph{3-interval exchange transformation} (3iet)
  with parameters $\alpha,\beta,\gamma$.
\end{defi}

With a 3-interval exchange transformation $T$, one can naturally
associate a ternary biinfinite word $u_T(x_0)=(u_n)_{n\in\Z}$,
which codes the orbit of a point $x_0$ from the domain of $T$, as
\begin{equation}\label{eq:u_T}
u_n = \begin{cases}
A & \text{if $\ T^n(x_0)\in I_A$,} \\
B & \text{if $\ T^n(x_0)\in I_B$,} \\
C & \text{if $\ T^n(x_0)\in I_C$.}
\end{cases}
\end{equation}

Similarly as in the case of a 2-interval exchange transformation,
the infinite word coding a 3iet can be periodic or aperiodic,
according to the choice of parameters $\alpha,\beta,\gamma$.
In paper~\cite{ambroz-masakova-pelantova-tcs}
it is proved that $u_T(x_0)$ is periodic if and only if $\alpha+\beta$ and
$\beta+\gamma$ are linearly dependent over $\Q$. Analogous to the
case of 2-interval exchange transformations, where periodic words
are not viewed as sturmian words, we will not consider periodic words.

\begin{defi}
  A biinfinite word $u_T(x_0)$ given by the prescription~\eqref{eq:u_T}
  is called a \emph{3iet word} with parameters $\alpha$,$\beta$,$\gamma$
  and $x_0$ if $\alpha+\beta$ and $\beta+\gamma$ are linearly independent over $\Q$.
\end{defi}

\subsection{Morphisms and incidence matrices}

A mapping $\varphi:\Acal^*\rightarrow\Acal^*$ is said to be a
\emph{morphism} over $\Acal$ if
$\varphi(w\widehat{w})=\varphi(w)\varphi(\widehat{w})$ holds for
any pair of finite words $w,\widehat{w}\in\Acal^*$. Obviously, a
morphism is uniquely determined by the images $\varphi(a)$ for all
letters $a\in\Acal$.
The action of a morphism $\varphi$ can be naturally extended to biinfinite words
by the prescription
\[
\varphi(u) = \varphi(\cdots u_{-2}u_{-1}u_0u_1\cdots) \coloneq
\cdots\varphi(u_{-2})\varphi(u_{-1})\varphi(u_{0})\varphi(u_{1})\cdots\,.
\]

To each morphism $\varphi$ over a $k$-letter alphabet
$\{a_1,\ldots,a_k\}$ one can assign its incidence matrix
$\Mcal_\varphi\in\N^{k\times k}$ by putting
\begin{equation}\label{eq:subst-matrix}
(\Mcal_\varphi)_{ij} = \text{number of letters $a_j$ in the word
$\varphi(a_i)$}\,.
\end{equation}
Morphisms over $\Acal$ form a monoid, whose neutral element is the
identity morphism. Let $\varphi$ and $\psi$ be morphisms over
$\Acal$, then the matrix of their composition, that is, of the
morphism $u\mapsto(\varphi\circ\psi)(u) =
\varphi\bigl(\psi(u)\bigr)$ is obtained by
\begin{equation}\label{eq:matrix-compose}
\mat{M}_{\varphi\circ\psi} = \mat{M}_{\psi}\mat{M}_{\varphi}\,.
\end{equation}
Therefore the mapping $\Rcal:\varphi\mapsto\mat{M}_\varphi$ is the
matrix representation of the monoid of all morphisms over $\Acal$.

\begin{defi}
  A morphism over the alphabet $\{A,B,C\}$ is said to be \emph{3iet preserving}
  if $\varphi(u)$ is a 3iet word for each 3iet word $u$. Monoid of all such
  morphisms will be denoted by $\PhiT$.
\end{defi}

\begin{ex}\label{ex:prehod-AC}
  It is easy to see that the morphism $\xi$ over $\{A,B,C\}$ given by prescriptions
  \begin{equation}\label{eq:pr-xi}
    A \mapsto C\,, \qquad B \mapsto B\,, \qquad C \mapsto A\,,
  \end{equation}
  is a 3iet preserving morphism. To a 3iet word, which codes the orbit of $x_0$
  under the transformation $T$ with intervals
  $[0,\alpha)\cup[\alpha,\alpha+\beta)\cup[\alpha+\beta,\alpha+\beta+\gamma)$,
  it assigns a 3iet word, which codes the orbit of $\alpha+\beta+\gamma-x_0$
  under the transformation $\tilde{T}$ with intervals
  $(0,\gamma]\cup(\gamma,\gamma+\beta]\cup(\gamma+\beta,\gamma+\beta+\alpha]$.
  The matrix of this morphism is
  $\mat{M}_\xi=\Big(\begin{smallmatrix}0&0&1\\0&1&0\\1&0&0\end{smallmatrix}\Big)$.
\end{ex}

In the proof of the fact that $\PhiT$ is not a finitely generated monoid we will
extensively use the following theorem, proved in~\cite{ambroz-masakova-pelantova-tcs}.

\begin{thm}
  Let $\mat{M}$ be a non-singular matrix of a 3iet preserving morphism. Then
  $\mat{M}$ belongs to the monoid
  $\EtN \coloneq \{\boldsymbol{M}\in\mathbb{N}^{3\times 3}\;|\;
  \boldsymbol{M}\boldsymbol{E}\boldsymbol{M}^T = \pm\boldsymbol{E}\:
  \text{ and }\det\boldsymbol{M}=\pm 1\}$, where $\boldsymbol{E} =
  \Big(\!\begin{smallmatrix}0&1&1\\-1&0&1\\-1&-1&0\end{smallmatrix}\Big)$.
\end{thm}

Let us denote
\[
\Rcal(\PhiT)\coloneq
\{\mat{M}_\varphi\:|\:\varphi\in\PhiT\text{ and }\det\mat{M}_\varphi\neq0\}\,.
\]
Obviously, $\Rcal(\PhiT)$ is a monoid too, and the above mentioned result can
be expressed as $\Rcal(\PhiT)\subset\EtN$. Moreover, the inclusion is strict,
that is, it is known that $\Rcal(\PhiT)\subsetneqq\EtN$.

We will prove that none of the monoids $\PhiT$, $\Rcal(\PhiT)$ or
$\EtN$ is finitely generated; for this we have to introduce the notion of
a non-decomposable element of a monoid.

\begin{defi}
  Let $\Mcal$ be a monoid with an operation $*$. An element $U\in\Mcal$ is called
  \emph{unit of $\Mcal$} if
  for any element $A\in\Mcal$ there exist $B,C\in\Mcal$ such that
  $A=B*U=U*C$. An element $M\in\Mcal$ is said to be
  \emph{non-decomposable in $\Mcal$}
  if $M=A*B$, $A,B\in\mathcal{M}$ implies that either $A$ or $B$ is unit of $\Mcal$.
\end{defi}

\section{Monoid $\EtN$ and its units}

At first we state the following necessary condition of a matrix to
be an element of $\EtN$.

\begin{prop}\label{prop:left_eivenvector}
  Let $\mat{M}\in\N^{3\times 3}$. If $\mat{M}\in\EtN$ then
  $(1,-1,1)$ is left eigenvector of $\mat{M}$.
\end{prop}
\begin{proof}
The matrix $\mat{E}$ has left eigenvector $(1,-1,1)$ with simple eigenvalue $0$.
As $\mat{M}$ is by definition of $\EtN$ non-singular, $(1,-1,1)$ is also left
eigenvector of $\mat{E}\mat{M}^T$ with simple eigenvalue 0. On the other hand
we have $(1,-1,1)\mat{M}\mat{E}\mat{M}^T=(1,-1,1)\mat{E}=(0,0,0)$, which means
that $(1,-1,1)\mat{M}$ is left eigenvector of $\mat{E}\mat{M}^T$ with simple
eigenvalue 0 as well. This implies $(1,-1,1)\mat{M}=c(1,-1,1)$.
\end{proof}

Let us describe in general units of a monoid $\Mcal$ of non-negative integral
matrices. By definition if $\mat{U}$ is a unit of $\Mcal$ then for the identity
matrix $\mat{I}\in\Mcal$ we have $\mat{I}=\mat{B}\mat{U}$, where $\mat{B}\in\Mcal$.
Hence the inverse matrix $\mat{U}^{-1}$ has to be also integral and non-negative.
This property holds only for permutation matrices; the only permutation
matrices fulfilling the necessary condition from
Proposition~\ref{prop:left_eivenvector} are
$\Big(\begin{smallmatrix}1&0&0\\0&1&0\\0&0&1\end{smallmatrix}\Big)$
and $\Big(\begin{smallmatrix}0&0&1\\0&1&0\\1&0&0\end{smallmatrix}\Big)$.
Recall that the second matrix is the matrix of the 3iet preserving
morphism $\xi$ from Example~\ref{ex:prehod-AC}, and, moreover,
$\xi$ is the only one morphism having this matrix. We can summarize these facts
in the following corollary.

\begin{coro}
\begin{enumerate}[(i)]
\item
  There are exactly two units in the monoids $\EtN$ and $\Rcal(\PhiT)$, namely the matrices
  $\Big(\begin{smallmatrix}1&0&0\\0&1&0\\0&0&1\end{smallmatrix}\Big)$
    and $\Big(\begin{smallmatrix}0&0&1\\0&1&0\\1&0&0\end{smallmatrix}\Big)$.
\item
  There are exactly two units in the monoid of 3iet preserving morphisms,
  namely the identity morphism and the morphism $\xi$.
\end{enumerate}
\end{coro}

To prove that the monoids $\EtN$, $\PhiT$ and $\Rcal(\PhiT)$ are not finitely
generated we will use the following Lemma.

\begin{lem}\label{lem:Mk_is_uniquely_decomp}
Let $k\geq 1$. Then matrix
$\mat{M}_k \coloneq \Big(\begin{smallmatrix}0&k&k-1\\0&k+1&k\\1&0&2\end{smallmatrix}\Big)$
is uniquely decomposable in $\EtN$, up to multiples of units.
\end{lem}
\begin{proof}
Let $k\geq 1$ be fixed, and let us consider $\mat{A},\mat{B}\in\EtN$,
$\mat{A},\mat{B}$ not units of $\EtN$, such that $\mat{M}_k=\mat{A}\mat{B}$.
We denote elements of $\mat{A}$, $\mat{B}$ and $\mat{M}_k$ by $a_{ij}$, $b_{ij}$
and $m_{ij}$, respectively.

From the null coefficients in the first column of $\mat{M}_k$ one
obtains $a_{11}b_{11}=a_{12}b_{21}=a_{13}b_{31}=a_{21}b_{11}=a_{22}b_{21}=a_{23}b_{31}=0$.
Let us assume that there is only one zero in the first column of $\mat{B}$, say $b_{i1}=0$.
This implies that $a_{1j}=a_{1k}= 0$ and $a_{2j}=a_{2k} = 0$ for $j,k\neq i$.
Hence $\mat{A}$ is singular, which is in contradiction with $\mat{A}\in\EtN$.
Therefore there are two zeros in $\mat{B}_{\bullet1}$. We denote by $\mat{B}_{\bullet i}$ and
$\mat{B}_{i\bullet}$ denote the $i$-th column and row of the matrix $\mat{B}$, respectively.
Recall that
$\mat{P}=\Big(\begin{smallmatrix}0&0&1\\0&1&0\\1&0&0\end{smallmatrix}\Big)=\mat{P}^{-1}$
is unit of $\EtN$, thus
it suffices to discuss only two cases for the first column of $\mat{B}$, namely
$b_{11}\neq 0$ and $b_{21}\neq 0$.

\medskip
\noindent\emph{Case A)} $b_{11}\neq 0$. We have
\[
\begin{pmatrix}0&a_{12}&a_{13}\\0&a_{22}&a_{23}\\a_{31}&a_{32}&a_{33}\end{pmatrix}
\begin{pmatrix}b_{11}&b_{12}&b_{13}\\0&b_{22}&b_{23}\\0&b_{32}&b_{33}\end{pmatrix}
=\begin{pmatrix}0&k&k-1\\0&k+1&k\\1&0&2\end{pmatrix}\,.
\]
Further elements of $\mat{A}$ and $\mat{B}$ can be enumerated by means of the following
facts \vspace{-0.4em}
\begin{itemize}
\item
  $m_{31} = 1 = a_{31}b_{11} \ \Rightarrow\ a_{31}=b_{11}=1$ \vspace{-0.4em}
\item
  $m_{32} = 0 = b_{12} + a_{32}b_{22} + a_{33}b_{32} \ \Rightarrow\ b_{12} = 0$ \vspace{-0.4em}
\item
  Since $\mat{B}\in\EtN$ the vector $(1,-1,1)$ has to be its left-eigenvector associated with
  an eigenvalue $\lambda$. Multiplying the first and the second column of $\mat{B}$ by
  this vector one has $(1,-1,1)\mat{B}_{\bullet1}=\lambda\ \Rightarrow\ \lambda=1$ and
  $(1,-1,1)\mat{B}_{\bullet2}=-\lambda\ \Rightarrow\ b_{22}=b_{32}+1$. \vspace{-0.4em}
\item
  Finally, $m_{32} = 0 = a_{32}(b_{32}+1)+a_{33}b_{32}\ \Rightarrow\ a_{32}=0$ and
  $a_{33}b_{32}=0$.
\end{itemize}
According to the term $a_{33}b_{32}=0$ we divide the discussion into two cases.

\medskip
\noindent\emph{Case A1)} $b_{32}=0$. We have
\[
\begin{pmatrix}0&a_{12}&a_{13}\\0&a_{22}&a_{23}\\1&0&a_{33}\end{pmatrix}
\begin{pmatrix}1&0&b_{13}\\0&1&b_{23}\\0&0&b_{33}\end{pmatrix}
=\begin{pmatrix}0&k&k-1\\0&k+1&k\\1&0&2\end{pmatrix}\,.
\]
As $\mat{B}$ is an upper triangular matrix, the product of its diagonal elements is equal
to its determinant, and we have $b_{33}=1$ and $1$ is the only eigenvalue of $\mat{B}$.
Further $m_{12}=k=a_{12}$ and
$m_{22}=k+1=a_{22}$. As before multiplying the third column of $\mat{B}$ by the left eigenvector
$(1,-1,1)$ associated to the eigenvalue $1$, one has
$b_{13}-b_{23}+1=1\ \Rightarrow\ b_{13}=b_{23}$. Thus $b_{23}>0$, otherwise
$\mat{B}$ is the identity matrix and so a unit of $\EtN$. From $m_{13}=k-1$ we get
$k-1=kb_{23}+a_{13}$, which is impossible for $b_{23}\geq 1$.

\medskip
\noindent\emph{Case A2)} $a_{33}=0$. We have
\[
\begin{pmatrix}0&a_{12}&a_{13}\\0&a_{22}&a_{23}\\1&0&0\end{pmatrix}
\begin{pmatrix}1&0&b_{13}\\0&b_{32}+1&b_{23}\\0&b_{32}&b_{33}\end{pmatrix}
=\begin{pmatrix}0&k&k-1\\0&k+1&k\\1&0&2\end{pmatrix}\,.
\]
Since $\mat{A}\in\EtN$ the vector $(1,-1,1)$ is its left-eigenvector.
Hence
$(1,-1,1)\mat{A}_{\bullet1}=\lambda\ \Rightarrow\ \lambda=1$,
$(1,-1,1)\mat{A}_{\bullet2}=-1\ \Rightarrow\ a_{22}=a_{12}+1$ and
$(1,-1,1)\mat{A}_{\bullet3}=1\ \Rightarrow\ a_{13}=a_{23}+1$. Therefore we have
$\pm1 = \det\mat{A} = a_{12}a_{23}-(a_{23}+1)(a_{12}+1) = a_{23}+a_{12}+1$. This is not
possible if $\det\mat{A}=-1$ and for $\det\mat{A}=1$ the equality can be fulfilled only
if $a_{12}=a_{23}=1$, i.e., $\mat{A} = \mat{P}$, which is unit of $\EtN$.

\medskip
\noindent\emph{Case B)} $b_{21}\neq 0$. We have
\[
\begin{pmatrix}a_{11}&0&a_{13}\\a_{21}&0&a_{23}\\a_{31}&a_{32}&a_{33}\end{pmatrix}
\begin{pmatrix}0&b_{12}&b_{13}\\b_{21}&b_{22}&b_{23}\\0&b_{32}&b_{33}\end{pmatrix}
=\begin{pmatrix}0&k&k-1\\0&k+1&k\\1&0&2\end{pmatrix}\,.
\]
Furthermore \vspace{-0.4em}
\begin{itemize}
\item
  $m_{31} = 1 = a_{32}b_{21}\ \Rightarrow\ a_{32}=b_{21}=1$ \vspace{-0.4em}
\item
  $m_{32} = 0 = a_{31}b_{12} + b_{22} + a_{33}b_{32}\ \Rightarrow\ b_{22}=0$ \vspace{-0.4em}
\item
  $(1,-1,1)\mat{B}_{\bullet1} = \lambda\ \Rightarrow\ \lambda = -1$ and
  $(1,-1,1)\mat{B}_{\bullet2} = 1\ \Rightarrow\ b_{12}+b_{32} = 1$. Without loss of
  generality $b_{12}=0$ and $b_{32}=1$. \vspace{-0.4em}
\item
  $m_{12} = k = a_{13}$, $m_{22}=k+1=a_{23}$ and $m_{32}=0=a_{33}$. \vspace{-0.4em}
\item
  $\pm1=\det\mat{B}=b_{13}\ \Rightarrow\ b_{13}=1$. \vspace{-0.4em}
\item
  $m_{33} = 2 = a_{31}+b_{23}$ and
  $(1,-1,1)\mat{B}_{\bullet3}=-1\ \Rightarrow\ b_{23}=b_{33}+2$ put together give
  $2 = a_{31} + b_{33} + 2\ \Rightarrow\ a_{31}=b_{33}=0$.
\end{itemize}
At this point the matrix $\mat{B}$ is fully enumerated, it is easy to
compute the remaining elements of $\mat{A}$, $a_{11}=k-1$ and $a_{21}=k$.

Therefore the matrix $\mat{M}_k$ has the following unique decomposition in $\EtN$
\begin{equation}\label{eq:decom_of_Mk}
\mat{M}_k = \begin{pmatrix}0&k&k-1\\0&k+1&k\\1&0&2\end{pmatrix} =
\underbrace{\begin{pmatrix}k-1&0&k\\k&0&k+1\\0&1&0\end{pmatrix}}_{\widetilde{\mat{M}}_k}
\underbrace{\begin{pmatrix}0&0&1\\1&0&2\\0&1&0\end{pmatrix}}_{\widetilde{\mat{M}}_1}
\,.
\end{equation}
\end{proof}

The uniqueness of the decomposition of $\mat{M}_k$ means that $\widetilde{\mat{M}}_k$
and $\widetilde{\mat{M}}_1$ in~(\ref{eq:decom_of_Mk}) are non-decomposable.
Hence one gets the following corollary.

\begin{coro}
The monoid $\EtN$ is not finitely generated.
\end{coro}

\section{Monoid of 3iet preserving morphisms}

The matrices $\mat{M}_k$ considered in Lemma~\ref{lem:Mk_is_uniquely_decomp} are
elements of $\Rcal(\PhiT)$ as shows the following lemma.

\begin{lem}\label{lem:phik_is_3iet}
  Let $k\in\N\setminus\{0\}$ and $\varphi_k:\{A,B,C\}^*\rightarrow\{A,B,C\}^*$ be
  the morphism given by
  \begin{equation}\label{eq:phik}
  A\mapsto B(CB)^{k-1},\quad B\mapsto B(CB)^k,\quad C\mapsto CAC\,.
  \end{equation}
  Then $\varphi_k$ is 3iet preserving.
\end{lem}
\begin{proof}
  Let us consider an arbitrary 3iet word $u$ with parameters $\alpha,\beta,\gamma$ and $x_0$.
  The corresponding 3-interval exchange transformation is given by
  \[
  T(x)=\begin{cases}
    x+\alpha+\beta & \text{if $x\in[0,\alpha)\eqcolon I_A$,}\\
    x+\alpha-\gamma & \text{if $x\in[\alpha,\alpha+\beta)\eqcolon I_B$,}\\
    x-\beta-\gamma & \text{if $x\in[\alpha+\beta,\alpha+\beta+\gamma)\eqcolon I_C$.}
  \end{cases}
  \]
  We show that the 3iet word with parameters $|I_A'|=\gamma$, $|I_B'|=k\alpha+(k+1)\beta$,
  $|I_C'|=(k-1)\alpha+k\beta+2\gamma$ and $x_0'=x_0$, where $I_A\cup I_B\subset I_B'$ and
  $I_C\subset I_C'$ (see Figure~\ref{fig:varphi_k}) coincides with the word $\varphi_k(u)$.

  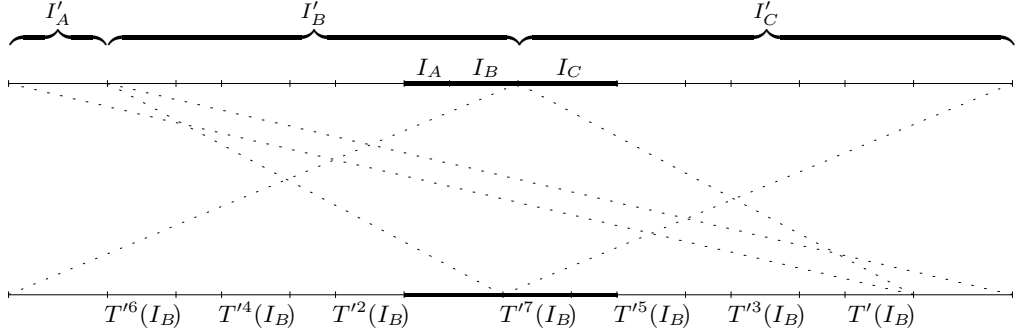
\begin{figure}[ht!]
    \centering
    \setlength{\unitlength}{0.4cm}
    \begin{picture}(35,10.9)(0,0)
      \drawline(1,8)(34,8)
        \drawline(1,7.9)(1,8.1) \drawline(4.25,7.9)(4.25,8.1)
        \drawline(6.5,7.9)(6.5,8.1) \drawline(8,7.9)(8,8.1)
        \drawline(10.25,7.9)(10.25,8.1) \drawline(11.75,7.9)(11.75,8.1)
        \drawline(14,7.9)(14,8.1) \drawline(15.5,7.9)(15.5,8.1)
        \drawline(17.75,7.9)(17.75,8.1)
        \drawline(21,7.9)(21,8.1) \drawline(23.25,7.9)(23.25,8.1)
        \drawline(24.75,7.9)(24.75,8.1) \drawline(27,7.9)(27,8.1)
        \drawline(28.5,7.9)(28.5,8.1) \drawline(30.75,7.9)(30.75,8.1)
        \drawline(34,7.9)(34,8.1)
        \drawline(1,1)(34,1)
        \drawline(1,0.9)(1,1.1) \drawline(4.25,0.9)(4.25,1.1)
        \drawline(6.5,0.9)(6.5,1.1) \drawline(8,0.9)(8,1.1)
        \drawline(10.25,0.9)(10.25,1.1) \drawline(11.75,0.9)(11.75,1.1)
        \drawline(14,0.9)(14,1.1) \drawline(17.25,0.9)(17.25,1.1)
        \drawline(19.5,0.9)(19.5,1.1)
        \drawline(21,0.9)(21,1.1) \drawline(23.25,0.9)(23.25,1.1)
        \drawline(24.75,0.9)(24.75,1.1) \drawline(27,0.9)(27,1.1)
        \drawline(28.5,0.9)(28.5,1.1) \drawline(30.75,0.9)(30.75,1.1)
        \drawline(34,0.9)(34,1.1)
       \dashline{0.1}(1,8)(30.75,1) \dashline{0.1}(4.25,8)(34,1)
       \dashline{0.1}(4.25,8)(17.25,1) \dashline{0.1}(17.75,8)(30.75,1)
       \dashline{0.1}(17.75,8)(1,1) \dashline{0.1}(34,8)(17.25,1)
       \put(1.05,9.2){$\overbrace{\hspace{1.26cm}}^{I_A'}$}
       \put(4.3,9.2){$\overbrace{\hspace{5.4cm}}^{I_B'}$}
       \put(17.8,9.2){$\overbrace{\hspace{6.5cm}}^{I_C'}$}
       \put(14.4,8.3){\text{\footnotesize$I_A$}}
       \put(16.25,8.3){\text{\footnotesize$I_B$}}
       \put(19,8.3){\text{\footnotesize$I_C$}}
       \put(28.6,0.1){\text{\footnotesize$T'(I_B\!)$}}
       \put(11.65,0.1){\text{\footnotesize$T'^2(I_B\!)$}}
       \put(24.65,0.1){\text{\footnotesize$T'^3(I_B\!)$}}
       \put(7.9,0.1){\text{\footnotesize$T'^4(I_B\!)$}}
       \put(20.9,0.1){\text{\footnotesize$T'^5(I_B\!)$}}
       \put(4.15,0.1){\text{\footnotesize$T'^6(I_B\!)$}}
       \put(17.15,0.1){\text{\footnotesize$T'^7(I_B\!)$}}
       \linethickness{0.4mm}
       \put(14,8){\line(1,0){7}}
       \put(14,1){\line(1,0){7}}
    \end{picture}
    \caption{The transformation $T'$ for $k=3$.}\label{fig:varphi_k}
  \end{figure}

  The transformation $T'$ corresponding to $|I_A'|$,$|I_B'|$ and $|I_C'|$ is the
  following one
  \[
  T'(x)=\begin{cases}
    x+|I_A'|+|I_B'| &
      \text{if $x\in[(1-k)\alpha-k\beta-\gamma,(1-k)\alpha-k\beta)\eqcolon I_A'$,}\\
    x+|I_A'|-|I_C'| & \text{if $x\in[(1-k)\alpha-k\beta,\alpha+\beta)\eqcolon I_B'$,}\\
    x-|I_B'|-|I_C'| & \text{if $x\in[\alpha+\beta,k\alpha+(k+1)\beta+2\gamma)\eqcolon I_C'$.}
  \end{cases}
  \]

  For a point $x\in I_C$ we have
  \begin{align*}
    T'(x) &= x - |I_A'| - |I_B'| =
        x - \gamma -k\alpha -(k+1)\beta \in I_A'\,,\\
    (T')^2(x) &= x + |I_C'| - |I_A'| =
        x + (k-1)\alpha + k\beta + \gamma \in I_C'\,,\\
    (T')^3(x) &= x + |I_C'| - 2|I_A'| - |I_B'| =
        x -\alpha -\beta = T(x)\,.
  \end{align*}
  Hence for any point $x\in I_C\subset I_C'$ its first iteration under $T'$ is
  $T'(x)\in I_A'$, the second iteration is $(T')^2(x)\in I_C'$ and the third one
  $(T')^3(x)$ sends it to the same place as the first iteration of the original transformation
  $T$. Therefore we substitute $C\mapsto CAC$.

  For a point $x\in I_B$ we successively have
  \begin{align*}
    T'(x) &= x + |I_C'| - |I_A'| = x + (k-1)\alpha + k\beta + \gamma \in I_C'\,, \\
    (T')^2(x) &= x + |I_C'| - 2|I_A'| - |I_B'| = x + -\alpha - \beta \in I_B'\,, \\
    (T')^3(x) &= x + 2|I_C'| - 3|I_A'| -|I_B'| = x + (k-2)\alpha + (k-1)\beta + \gamma \in I_C'\,, \\
    (T')^4(x) &= x + 2|I_C'| - 4|I_A'| - 2|I_B'| = x + -2\alpha - 2\beta \in I_B'\,,
  \end{align*}
  that is, in general
  \begin{align*}
    (T')^{2l+1}(x) &= x + (l+1)|I_C'| - (2l+1)|I_A'| - l|I_B'| = \\
                   &= x + (k-l-1)\alpha + (k-l)\beta + \gamma \in I_C'\,, \\
    (T')^{2l}(x) &= x + l|I_C'| - 2l|I_A'| - l|I_B'| = x + -l\alpha - l\beta \in I_B'\,.
  \end{align*}
  The iterations belong alternately to intervals $I_C'$ and $I_B'$, with the last one being
  $(T')^{2k+1}(x)=x-\alpha+\gamma=T(x)$. Therefore we substitute $B\mapsto B(CB)^k$.

  For a point $x\in I_A$ the iterations are the same as in the case of $x\in I_B$ with the
  last one being $(T')^{2(k-1)+1}(x) = x + \beta + \gamma = T(x)$, and hence
  $A\mapsto B(CB)^{k-1}$.
  It follows from the definition that $\mat{M}_k$ is the matrix of morphism
  $\varphi_k$ given by these three assignments.
\end{proof}

\begin{thm}
  The monoid of 3iet preserving morphisms $\PhiT$ is not finitely generated.
\end{thm}
\begin{proof}
  At first note that there is an obvious relation between the decomposition of a
  morphism and of its matrix. Let $\varphi\in\PhiT$ be a morphism having non-trivial
  decomposition $\varphi=\psi\circ\phi$. Then its matrix $\mat{M}_\varphi$ is
  decomposable in $\Rcal(\PhiT)$, $\mat{M}_\varphi = \mat{M}_\phi\mat{M}_\psi$
  with $\mat{M}_\phi,\mat{M}_\psi$ not being units of $\Rcal(\PhiT)$.

  To prove that $\PhiT$ is not finitely generated we construct
  the sequence $(\widetilde{\varphi}_k)_{k\in\N}$ of non-decomposable
  3iet preserving morphisms. To meet this aim we exploit the morphism
  $\varphi_k\in\Phi_{3\mathrm{iet}}$, defined in
  Lemma~\ref{lem:phik_is_3iet}. There are two possible cases, depending on the
  nature of the matrices $\widetilde{\mat{M}}_1$ and $\widetilde{\mat{M}}_k$, to
  which $\mat{M}_k$ (i.e., the matrix of $\varphi$) can be decomposed
  (cf.\ proof of Lemma~\ref{lem:Mk_is_uniquely_decomp})
  \begin{enumerate}[(i)]
  \item
    If $\widetilde{\mat{M}}_1$ is not the matrix of a 3iet preserving morphism
    then $\varphi_k$ are non-decomposable morphisms for all $k\in\N$, and therefore
    if suffices to put $\widetilde{\varphi}_k = \varphi_k$.
  \item
    Let us suppose that $\widetilde{\mat{M}}_1$ is the matrix of a 3iet preserving
    morphism. We define the sequence $(\widetilde{\varphi}_k)_{k\in\N}$ in the
    following way. If $\widetilde{\mat{M}}_k$ is the matrix of a 3iet preserving
    morphism $\eta_k$ then $\widetilde{\varphi}_k = \eta_k$, otherwise
    $\widetilde{\varphi}_k = \varphi_k$.
  \end{enumerate}
\end{proof}

\begin{coro}
The monoid $\Rcal(\PhiT)$ is not finitely generated.
\end{coro}


\section*{Acknowledgements}

The authors acknowledge financial support by Czech Science
Foundation GA \v{C}R 201/05/0169, and by the grant LC06002 of the
Ministry of Education, Youth, and Sports of the Czech Republic.



\end{document}